# Study of Some Improved Ratio Type Estimators Using Information on Auxiliary Attributes Under Second Order Approximation

*Prayas Sharma, *Rajesh Singh and **Jong-Min Kim

*Department of Statistics, Banaras Hindu University, Varanasi(U.P.)-221005, India

**Statistics, Division of Science and Mathematics, University of Minnesota- Morris

(prayassharma02@gmail.com, rsinghstat@gmail.com, jongmink@morris.umn.edu )

Abstract

Chakrabarty (1979), Khoshnevisan et al. (2007), Sahai and Ray (1980), Solanki et al. (2012) suggested some estimators to estimate unknown population mean of the study variable. These authors discussed the estimators along with their first order biases and mean square errors(MSE's). In this paper, we have tried to found out the second order biases and mean square errors of some estimators using information on auxiliary attribute. We have compared the performance of the estimators with the help of a numerical illustration.

## 1. INTRODUCTION

In the theory of sample surveys, use of the auxiliary information can increase the precision or accuracy of an estimator of unknown population parameter of interest when study variable y is highly correlated with the auxiliary variable x. But there may be many practical situations when auxiliary information is not available directly (it is qualitative in nature), that is, auxiliary information is available in the form of an attribute. For example:

(a) The height of a person may depend on the fact that whether the person is male or female.

(b) The efficiency of a Dog may depend on the particular breed of that Dog.

(c) The yield of wheat crop produced may depend on a particular variety of wheat, etc.

In these situations by taking the advantage of point bi-serial correlation between the study variable y and the auxiliary attribute $\phi$ along with the prior knowledge of the population parameter of auxiliary attribute, the estimators of population parameter of interest can be constructed.

Consider a sample of size n drawn by simple random sampling without replacement (SRSWOR) from a population of size N. let $y_i$ and $\phi_i$ denote the observation on variable y and x respectively for the $i^{th}$ unit (i=1,2,3,……N) . We note that $\phi_i$ =1, if $i^{th}$ unit possesses attributes $\phi$ and $\phi_i$ =0 otherwise . Let $A = \sum_{i=1}^{N} \phi_i$, and $a = \sum_{i=1}^{n} \phi_i$ denote the total no. of units in the population and sample respectively possessing attribute x. Let $P = \frac{A}{N}$ and $p = \frac{a}{n}$ denote the proportion of units in the population and sample respectively possessing attribute $\phi$.

Using the information of point biserial correlation between the study variable and the auxiliary attribute Naik and Gupta (1996), Shabbir and Gupta (2006), Ab-Alfatah *et al.* (2010) and Singh *et al.* (2007, 2008 ) have suggested improved estimators for estimating unknown population mean $\bar{Y}$. In this paper we have studied properties of some estimators under second order of approximation.

## 2. Some Estimators in Simple Random Sampling

For estimating the population mean $\bar{Y}$ of Y, adapting Chakrabarty (1979) ratio-type estimator in case of attribute auxiliary variable, we have

$$t_1 = (1-\alpha)\bar{y} + \alpha\bar{y}\frac{P}{p} \tag{2.1}$$

where $\bar{y} = \frac{1}{n}\sum_{i=1}^{n} y_i$ .

Khoshnevisan et al. (2007) ratio- type estimator in case of attribute auxiliary attribute is given by

$$t_2 = \bar{y}\left[\frac{P}{\beta p + (1-\beta)P}\right]^g \qquad (2.2)$$

where g is a constant. For g=1 , $t_2$ is same as conventional ratio estimator whereas for g = -1, it becomes conventional product type estimator.

Adapting Sahai and Ray (1980) estimator for the situation when information is available in form of auxiliary attribute, we have

$$t_3 = \bar{y}\left[2 - \left\{\frac{p}{P}\right\}^w\right] \qquad (2.3)$$

where w is a constant.

Adapting Solanki et al. (2012) estimator for case of attributes, we get an estimator $t_4$ as

$$t_4 = \bar{y}\left[2 - \left\{\left(\frac{p}{P}\right)^{\lambda} \exp\frac{\delta(p-P)}{(p+P)}\right\}\right] \qquad (2.5)$$

where $\lambda$ is a constant, suitably chosen by minimizing mean square error of the estimator $t_4$ .

**3. Notations used**

Let us define,

$$e_0 = \frac{\bar{y} - \bar{Y}}{\bar{Y}} \quad \text{and} \quad e_1 = \frac{p - P}{P},$$

then $E(e_0) = E(e_1) = 0$.

For obtaining the bias and MSE expressions of the estimators, following lemmas will be used:

**Lemma 3.1**

(i) $V(e_0) = E\{(e_0)^2\} = \frac{N-n}{N-1}\frac{1}{n}C_{02} = L_1C_{02}$

(ii) $V(e_1) = E\{(e_1)^2\} = \frac{N-n}{N-1}\frac{1}{n}C_{20} = L_1C_{20}$

(iii) $COV(e_0, e_1) = E\{(e_0e_1)\} = \frac{N-n}{N-1}\frac{1}{n}C_{11} = L_1C_{11}$

**Lemma 3.2**

(iv) $E\{(e_1{}^2e_0)\} = \frac{(N-n)}{(N-1)}\frac{(N-2n)}{(N-2)}\frac{1}{n^2}C_{21} = L_2C_{21}$

(v) $E\{(e_1{}^3)\} = \frac{(N-n)}{(N-1)}\frac{(N-2n)}{(N-2)}\frac{1}{n^2}C_{30} = L_2C_{30}$

**Lemma 3.3**

(vi) $E(e_1{}^3e_0\ ) = L_3C_{31} + 3L_4C_{20}C_{11}$

(vii) $E\{(e_1{}^4)\} = \frac{(N-n))(N^2+N-6nN+6n^2)}{(N-1)(N-2)(N-3)}\frac{1}{n^3}C_{30} = L_3C_{40} + 3L_4C_{20}{}^2$

(viii) $E(e_1^2e_0{}^2) = L_3C_{40} + 3L_4C_{20}$

Where $L_3 = \frac{(N-n))(N^2+N-6nN+6n^2)}{(N-1)(N-2)(N-3)}\frac{1}{n^3}$, $L_4 = \frac{N(N-n))(N-n-1)(n-1)}{(N-1)(N-2)(N-3)}\frac{1}{n^3}$

And $Cpq = \frac{(Xi - \bar{X})^p}{\bar{X}^p}\frac{(Y_i - \bar{Y})^q}{\bar{Y}^q}$

Proof of these lemma are straight forward by using SRSWOR ( see Sukhatme and Sukhatme (1970)).

## 4. First Order Biases and Mean Squared Errors

The expressions for the biases of the estimators $t_1$, $t_2$ $t_3$ and $t_4$ are, respectively written as

$$\text{Bias}(t_1) = \bar{Y}\left[\frac{1}{2}\alpha L_1C_{20} - \alpha L_1C_{11}\right] \tag{4.1}$$

$$\text{Bias}(t_2) = \overline{Y}\left[\frac{g(g+1)}{2}L_1C_{20} - g\beta L_1C_{11}\right] \quad (4.2)$$

$$\text{Bias}(t_3) = \overline{Y}\left[-\frac{w(w-1)}{2}L_1C_{20} - wL_1C_{11}\right] \quad (4.3)$$

$$\text{Bias}(t_4) = \overline{Y}\left[-\frac{K(K-1)}{2}L_1C_{20} - KL_1C_{11}\right] \quad (4.5)$$

where, $k = \frac{(\delta + 2\lambda)}{2}$.

The expressions for the MSE's of the estimators $t_{1,}$ $t_2$ $t_3$ and $t_4$ are , respectively, given by

$$\text{MSE}(t_1) = \overline{Y}^2\left[L_1C_{02} + \alpha^2 L_1C_{20} - 2\alpha L_1C_{11}\right] \quad (4.6)$$

$$\text{MSE}(t_2) = \overline{Y}^2\left[L_1C_{02} + g^2\beta^2 L_1C_{20} - 2g\beta L_1C_{11}\right] \quad (4.7)$$

$$\text{MSE}(t_3) = \overline{Y}^2\left[L_1C_{02} + w^2 L_1C_{20} - 2wL_1C_{11}\right] \quad (4.8)$$

$$\text{MSE}(t_4) = \overline{Y}^2\left[L_1C_{02} + k^2 L_1C_{20} - 2kL_1C_{11}\right] \quad (4.10)$$

The MSE's of the estimators $t_{1,}$ $t_2$ $t_3$ and $t_4$ under optimum conditions are equal to the MSE of the regression estimator when we consider the terms up to first order of approximations. In search of the optimum estimator, we have extended this study to second order of approximation.

## 5. Second Order Biases and Mean Squared Errors

Expressing estimator $t_1$ in terms of e's (i=0,1), we get

$$t_1 = \overline{Y}(1+e_0)\left\{(1-\alpha) + \alpha(1+e_1)^{-1}\right\}$$

Or

$$t_1 - \overline{Y} = \overline{Y}\left\{e_0 + \frac{e_1}{2} + \frac{\alpha}{2}e_1^{\ 2} - \alpha e_0 e_1 + \alpha e_0 e_1^{\ 2} - \frac{\alpha}{6}e^3 - \frac{\alpha}{6}e_0 e_1^{\ 3} + \frac{\alpha}{24}e^4\right\} \quad (5.1)$$

Taking expectations, we get the bias of the estimator $t_1$ up to the second order of approximation as

$$\text{Bias}_2(t_1) = \overline{Y}\left[\frac{\alpha}{2}L_1C_{20} - \alpha L_1C_{11} - \frac{\alpha}{6}L_2C_{30} + \alpha L_2C_{21} - \frac{\alpha}{6}\left(L_3C_{31} + 3L_4C_{20}C_{11}\right)\right.$$

$$\left. + \frac{\alpha}{24}(L_3C_{40} + 3L_4C_{20}{}^2)\right] \tag{5.2}$$

Similarly, we get the expressions for the biases of the estimator's $t_2$, $t_3$, and $t_4$ up to second order of approximation as follows

$$\text{Bias}_2(t_2) = \overline{Y}\left[\frac{g(g+1)}{2}\beta^2 L_1C_{20} - g\beta L_1C_{11} - \frac{g(g+1)}{2}\beta^2 L_2C_{21} - \frac{g(g+1)(g+2)}{6}\beta^3 L_2C_{30}\right.$$

$$- \frac{g(g+1)(g+2)}{6}\beta^3\left(L_3C_{31} + 3L_4C_{20}C_{11}\right)$$

$$\left. + \frac{g(g+1)(g+2)(g+3)}{24}\beta^4(L_3C_{40} + 3L_4C_{20}{}^2)\right] \tag{5.3}$$

$$\text{Bias}_2(t_3) = \overline{Y}\left[\frac{w(w-1)}{2}L_1C_{20} - wL_1C_{11} - \frac{w(w-1)}{2}L_2C_{21} - \frac{w(w-1)(w-2)}{6}L_2C_{30}\right.$$

$$- \frac{w(w-1)(w-2)}{6}\left(L_3C_{31} + 3L_4C_{20}C_{11}\right)$$

$$\left. - \frac{w(w-1)(w-2)(w-3)}{24}(L_3C_{40} + 3L_4C_{20}{}^2)\right] \tag{5.4}$$

$$\text{Bias}_2(t_4) = \overline{Y}\left[-\frac{k(k-1)}{2}L_1C_{20} - kL_1C_{11} - \frac{k(k-1)}{2}L_2C_{21} - ML_2C_{30} - M\left(L_3C_{31} + 3L_4C_{20}C_{11}\right)\right.$$

$$\left. - N(L_3C_{40} + 3L_4C_{20}{}^2)\right] \tag{5.6}$$

where, $M = \frac{1}{2}\left\{\frac{\left(\delta^3 - 6\delta^2\right)}{24} + \frac{\left(\alpha(\delta^2 - 2\delta\right)}{4} + \frac{\lambda(\lambda-1)}{2}\delta + \frac{\lambda(\lambda-1)(\lambda-2)}{3}\right\}$ and $k = \frac{(\delta + 2\lambda)}{2}$

$$N=\frac{1}{8}\left\{\frac{\left(\delta^4-12\delta^3+12\delta^2\right)}{48}+\frac{\left(\alpha(\delta^3-6\delta)\right)}{6}+\frac{\lambda(\lambda-1)}{2}(\delta^2-2\delta)+\frac{\lambda(\lambda-1)(\lambda-2)(\lambda-3)}{3}\right\}$$

Following are the expressions of the MSE's of the estimator's $t_1$, $t_2$, $t_3$ and $t_4$ respectively, up to second order of approximation

$$\begin{aligned}MSE_2(t_1)=\overline{Y}^2\Big[&L_1C_{02}+\alpha^2L_1C_{20}-2\alpha L_1C_{11}-\alpha^2L_2C_{30}+(2\alpha^2+\alpha)L_2C_{21}\\&-2\alpha^2\left(L_3C_{31}+3L_4C_{20}C_{11}\right)\\&+\alpha(\alpha+1)\left(L_3C_{22}+3L_4(C_{20}C_{02}+C_{11}^2)\right)+\frac{5}{24}\alpha^2(L_3C_{40}+3L_4C_{20}{}^2)\Big]\end{aligned} \quad (5.7)$$

$$\begin{aligned}MSE_2(t_2)=\overline{Y}^2\Big[&L_1C_{02}+g^2\beta^2L_1C_{20}-2\beta gL_1C_{11}-\beta^3g^2(g+1)L_2C_{30}+g(3g+1)\beta^2L_2C_{21}\\&-2\beta gL_2C_{12}-\left\{\frac{7g^3+9g^2+2g}{3}\right\}\beta^3\left(L_3C_{31}+3L_4C_{20}C_{11}\right)\\&+g(2g+1)\beta^2\left(L_3C_{22}+3L_4(C_{20}C_{02}+C_{11}^2)\right)\\&+\left\{\frac{2g^3+9g^2+10g+3}{6}\right\}\beta^4(L_3C_{40}+3L_4C_{20}{}^2\Big]\end{aligned} \quad (5.8)$$

$$\begin{aligned}MSE_2(t_3)=\overline{Y}^2\Big[&L_1C_{02}+w^2L_1C_{20}-2wL_1C_{11}-w^2(w-1)L_2C_{30}+w(w+1)L_2C_{21}-2wL_2C_{12}\\&+\left\{\frac{5w^3-3w^2-2w}{3}\right\}\left(L_3C_{31}+3L_4C_{20}C_{11}\right)\\&+w\left(L_3C_{22}+3L_4(C_{20}C_{02}+C_{11}^2)\right)+\left\{\frac{7w^4-18w^3+11w^2}{24}\right\}(L_3C_{40}+3L_4C_{20}{}^2)\Big]\end{aligned} \quad (5.9)$$

$$MSE_2(t_4)=\overline{Y}^2\Big[L_1C_{02}+k^2L_1C_{20}-2kL_1C_{11}+kL_2C_{21}-2kL_2C_{12}+k^2(k-1)L_2C_{30}$$

$$+2k^2(k-1)\left(L_3C_{31}+3L_4C_{20}C_{11}\right)+k\left(L_3C_{22}+3L_4(C_{20}C_{02}+C_{11}^2)\right)$$
$$\left.+\frac{(k^2-k)^2}{4}(L_3C_{40}+3L_4C_{20}{}^2)\right] \qquad (5.11)$$

where, $k=\frac{(\delta+2\lambda)}{2}$.

The optimum value of w for the constant involved in the estimator $t_3$, we get by minimizing $MSE_2(t_3)$. But theoretically the determination of the optimum value of w is very difficult, so we have calculated the optimum value by using numerical techniques. Similarly, the optimum value of k which minimizes the MSE of the estimator $t_4$ is obtained by using numerical techniques.

## 6. Numerical Illustration

The various result obtained in the pervious sections are now examined with the help of the following data set:

**Source of the Data**

The data for the empirical analysis is taken from Sukhatme and Sukhatme (1970), p.256

Y= number of villages in the circles and

$\phi$ = A circle consisting more than five villages

N=89, $\overline{Y}=3.36$, $P=0.1236$, $\rho_{pb}=0.766$, $C_y=0.604$, $C_x=2.19$.

**Table 6.1: Biases and MSE's of estimators**

| **Estimators** | **Bias** | | **MSE** | |
|---|---|---|---|---|
| | First order | Second order | First order | Second order |
| $t_1$ | 23.385388 | 24.42412947 | 17597.0515 | 23161.50715 |
| $t_2$ | -9.11939E-09 | -0.84723967 | 17597.0515 | 17653.75783 |
| $t_3$ | 18.30588379 | 26.44381147 | 17597.0515 | 19089.85861 |
| $t_4$ | 18.30588379 | 55.4734888 | 17597.0515 | 18104.07826 |

In the Table 6.1 the biases and MSE's of the estimator's $t_1$, $t_2$, $t_{3,}$ and $t_4$ are written under first order and second order of approximations. For all the estimators $t_1$, $t_2$, $t_{3,}$ and $t_4$ it is observed that the value of the biases are increasing. For all estimators MSE's up to the first order of approximation under optimum conditions are equal, which prompt us to study the properties of the estimators up to the second order of approximation. On the basis of study up to the second order of approximation we conclude that estimator $t_2$ is best followed by $t_4$, and $t_3$ among the estimators considered here for the given data set.

**REFERENCES**

Abd-Elfattah, A.M. El-Sherpieny, E.A. Mohamed, S.M. Abdou, O. F. (2010): Improvement in estimating the population mean in simple random sampling using information on auxiliary attribute. Appl. Math. and Compt. doi:10.1016/j.amc.2009.12.041

Chakrabarty, R.P. (1979) : Some ratio estimators. Journal of the Indian Society of Agricultural Statistics 31(1), 49–57.

Khoshnevisan, M., Singh, R., Chauhan, P., Sawan, N., and Smarandache, F. (2007). A general family of estimators for estimating population mean using known value of some population parameter(s), Far East Journal of Theoretical Statistics 22 181–191.

Naik, V.D and Gupta, P.C. (1996): A note on estimation of mean with known population proportion of an auxiliary character. Jour. Ind. Soc. Agri. Stat., 48(2), 151-158.

Ray, S.K. and Sahai, A (1980) : Efficient families of ratio and product type estimators, Biometrika 67(1), 211–215.

Shabbir, J. and Gupta, S. (2006): A new estimator of population mean in stratified sampling, Commun. Stat. Theo. Meth. 35: 1201–1209.

Singh, R., Cauhan, P., Sawan, N., and Smarandache, F. (2007). Auxiliary Information and A Priori Values in Construction of Improved Estimators. Renaissance High Press.

Singh, R., Chauhan, P., Sawan, N. and Smarandache,F. (2008): Ratio estimators in simple random sampling using information on auxiliary attribute. Pak. J. Stat. Oper. Res.,4,1,47-53.

Solanki, R.S., Singh, H. P. and Rathour, A. (2012) : An alternative estimator for estimating the finite population mean using auxiliary information in sample surveys. ISRN Probability and Statistics doi:10.5402/2012/657682

Sukhatme, P.V. and Sukhatme, B.V. (1970): Sampling theory of surveys with applications. Iowa State University Press, Ames, U.S.A.